\documentclass[12pt,a4paper]{article}
\usepackage[utf8]{inputenc}
\usepackage[english]{babel}
\usepackage{amsmath}
\usepackage{amsfonts}
\usepackage{amssymb}
\usepackage{makeidx}
\usepackage{graphicx}
\usepackage[left=2cm,right=2cm,top=2cm,bottom=2cm]{geometry}
\usepackage{amsthm}
\usepackage{xcolor}
\usepackage[all]{xy}
\usepackage{tikz}

\providecommand{\keywords}[1]
{
  \small	
  \textbf{Keywords:} #1 pro-$p$ groups, amalgamated free pro-$p$ products, pro-$p$ HNN-extensions
}

\title{Cyclic splittings of pro-$p$  groups}

\author{Jesus Berdugo $^{2}$, Pavel Zalesskii$^{1}$  \\
        \small $^{1}$ University of Brasilia, Department of mathematics,  E-mail: pz@mat.unb.br  \\
        \small $^{2}$ University of Brasilia, Department of mathematics,   E-mail: berdugo@mat.unb.br \\
}

\date{\today}

\newtheorem{teo}{Theorem}[section]

\newtheorem{theorem}[teo]{Theorem}
\newtheorem{proposition}[teo]{Proposition}%
\newtheorem{remark}[teo]{Remark}%
\newtheorem{corollary}[teo]{Corollary}

\newtheorem{definition}[teo]{Definition}%
\newtheorem{lemma}[teo]{Lemma}

\newcommand{\Z}{\mathbb{Z}}

\begin{document}
\maketitle

\abstract{In this paper, we prove a pro-$p$ version of the Rips-Sela's Theorems on  splittings of a group as an amalgamated free product or HNN-extension over an infinite cyclic subgroup.}

\vspace{0.5cm}


\section{Introduction}\label{sec1}

In  1997 Rips and Sela published the fundamental paper \cite{bib1}, where they studied infinite cyclic splittings (i.e. $\Z$-splittings) of groups as an amalgamated free product or an HNN-extension.  They  constructed a canonical JSJ decomposition  for finitely presented groups  that gives a complete description of all $\Z$-splittings of these groups.
In order to understand all possible $\Z$-splittings of a group, they needed to study carefully the "interaction" between any two given elementary $\Z$-splittings of it.
 
The objective of this paper is to study the "interaction" between any two given  $\Z_p$-splittings of a pro-$p$ group. Namely, we prove a pro-$p$ version of Rips-Sela's theorems on $\mathbb{Z}$-splittings (  \cite[Theorem 2.1 and Theorem 3.6]{bib1}).
 
A splitting of a pro-$p$ group $G$ as an amalgamated free pro-$p$ product or HNN-extension over $\mathbb{Z}_p$ will be called a $\mathbb{Z}_p$-splitting in the paper. An element of $G$ is called elliptic with respect to a splitting $G=G_1\amalg_{\Z_p} G_2$ as an amalgamated free pro-$p$ product (resp. pro-$p$ HNN-extension $G=HNN(G_1, \Z_p, t)$) if it is conjugate into $G_1\cup G_2$  (resp. into $G_1$) and is called hyperbolic otherwise. A pair of given $\mathbb{Z}_{p}$-splittings $A_1\amalg_{C_1} B_1$ and $A_2\amalg_{C_2} B_2$ over $C_1=\langle c_1\rangle$, $C_2=\langle c_2\rangle$ is called:

\begin{itemize}

\item $Elliptic-Elliptic:$ If $c_{1}$ is elliptic in $A_2\amalg_{C_2} B_2$ and $c_{2}$ is elliptic in $A_1\amalg_{C_1} B_1$.

\item $Hyperbolic-Hyperbolic:$ If $c_{1}$ is hyperbolic in $A_2\amalg_{C_2} B_2$ and $c_{2}$ is hyperbolic in $A_1\amalg_{C_1} B_1$.

\item $Hyperbolic-Elliptic:$ If $c_{1}$ is hyperbolic in $A_2\amalg_{C_2} B_2$ and $c_{2}$ is elliptic in $A_1\amalg_{C_1} B_1$.
\end{itemize}

\begin{definition} A pro-$p$ group $G$ is said to be \textit{freely indecomposable group} if it can not be written as a free pro-$p$ product of two non-trivial subgroups.
\end{definition}

Our first result is the pro-$p$ analog of \cite[Theorem 2.1]{bib1}.

\begin{theorem}\label{sela} Let $G$ be a finitely generated freely indecomposable pro-$p$ group. Then any two $\mathbb{Z}_{p}$-splittings of $G$ are either elliptic-elliptic or hyperbolic-hyperbolic.
\end{theorem}

 Next, given two hyperbolic-hyperbolic splittings over $C_1$ and $C_2$, we study the normalizer $N_G(C_i)$, $i=1,2$.  This study corresponds to the results of Section 3 of Rips-Sela's paper (where normalizers are called anti-centralizers). Note that Sela and Rips \cite{bib1} use  the existence of Tits' axis, on which the hyperbolic element  acts. In the pro-$p$ case, such an axis does not exist, so our argument is different from their argument. We divide our argument  into two cases: $p>2$ and $p=2$ since for $p=2$ the study
 of all possible normalizers is quite detailed and involves a sequence of case
 studies. We state here the case $p>2$, since it has only two cases, and refer the reader to Proposition \ref{prop1 p=2} for the $p=2$ case.
 
 \begin{proposition} \label{prop1 p>2} Let $p>2$ and suppose that $G=A_{1}\amalg_{C_{1}}B_{1}$ (or $G=HNN(A_{1},C_{1},t_{1})$) and $G=A_{2}\amalg_{C_{2}}B_{2}$ (or $G=HNN(A_{2},C_{2},t_{2})$)  are two hyperbolic-hyperbolic $\mathbb{Z}_{p}$-splittings.  Let $H_i\neq 1$ be a subgroup of $C_i$. Then  $N_{G}(H_{i})$  has one of the following types:

\begin{enumerate}
\item[(i)]  cyclic group $\Z_p$; 
\item[(ii)]  
$\Z_p\times\Z_p$; 

\end{enumerate}
\end{proposition}
 
Next we show that if $N_G(C_1)$ is not an infinite cyclic or infinite dihedral group and G does not split over a subgroup of order $\leq 2$, then $N_G(C_1)$ is in fact the ambient group $G$; hence,  for $p>2$ $G$ is $\Z_p\times \Z_p$ and for $p=2$ $G$ is $\Z_2\times \Z_2$ by $\Z/2\Z$. We state again here the $p>2$ case; for the $p=2$ case the reader can consult  Theorem \ref{teo8}.
 
\begin{theorem}\label{sela2} Let $p>2$ and
let  $G$ be a finitely generated  pro-$p$ group that does not split as a free pro-$p$ product.  Let $G=A_{1}\amalg_{C_{1}} B_{1}$ (or $G=HNN(A_{1},C_{1},t_{1})$), and $G=A_{2}\amalg_{C_{2}} B_{2}$ (or $G=HNN(A_{2},C_{2},t_{2})$) be two hyperbolic-hyperbolic $\mathbb{Z}_{p}$-splittings of $G$. Suppose that $N_{G}(C_1)$ is not cyclic. Then $G\cong \Z_p\times \Z_p$.
\end{theorem}
 
 We deduce the following 
 
\begin{corollary}\label{malnormal} With the hypothesis of Theortem \ref{sela2} suppose $G$ is non-abelian. Then $C_i$ is malnormal in either $A_i$ or $B_i$ (Resp. $C_i$ or $C_i^{t_i}$ is  malnormal  in $A_i$), $i=1,2$. 
\end{corollary}
 
As in classical Bass-Serre's theory, a pro-$p$ group that splits as an amalgamated free pro-$p$ product or HNN-extension acts on the corresponding pro-$p$ tree. The action is called $k$-acylindrical if any non-trivial element can fix at most $k$ consecutive edges. We deduce from Corollary \ref{malnormal} the $2$-acylindricity of the action. 
 
\begin{theorem}\label{2-acylindrical p>2}
Let  $p>2$ and $G$ be a non-abelian finitely generated  pro-$p$ group that does not split as a free pro-$p$ product.  Let $G=A_{1}\amalg_{C_{1}} B_{1}$ (or $G=HNN(A_{1},C_{1},t_{1})$), and $G=A_{2}\amalg_{C_{2}} B_{2}$ (or $G=HNN(A_{2},C_{2},t_{2})$) be two hyperbolic-hyperbolic $\mathbb{Z}_{p}$-splittings of $G$.  Then the action of $G$ on the standard pro-$p$ trees of these splittings is 2-acylindrical.
\end{theorem}

The proofs of the results use  the pro-$p$ version of the Bass-Serre theory  that can be found in \cite{bib5}.\\ 
 
The structure of the paper is as follows: In Section \ref{sec2} we give basic definitions that will be used throughout the article. In Section \ref{sec4} we prove Theorem \ref{sela} and  in  Section \ref{sec5} we prove Theorem \ref{sela2}. 
 
 \bigskip
 {\bf Conventions.} Throughout the paper, unless otherwise stated, groups are pro-$p$, subgroups are closed, and homomorphisms are continuous. In particular, $\langle S\rangle$ will mean the topological generation in the paper and presentations are taking in the category of pro-$p$ groups; $a^g$ will stand for $g^{-1}ag$ in the paper.

\section{Preliminaries}\label{sec2}
\subsection{Basic definitions }

\begin{definition}
A $graph$ $\Gamma$ is a disjoint union $E(\Gamma)\cup V(\Gamma)$, with the two maps $d_0,d_1:\Gamma\longrightarrow V(\Gamma)$, whose restriction to $V(\Gamma)$ are the identity map and for any element $e\in E(\Gamma)$, $d_0(e)$ and $d_1(e)$ are the initial and the terminal vertices of $e$ respectively. 
\end{definition}

\begin{definition}\label{prograph} A graph
$\Gamma$ is called a profinite graph if $\Gamma$ is a profinite space with a non-empty subset $V(\Gamma)$ such that:

\begin{itemize}
    \item[1.] $V(\Gamma)$ is closed,
    \item[2.] the maps $d_0,d_1:\Gamma\longrightarrow V(\Gamma)$ are continuous.
\end{itemize}
\end{definition}

We call $V (\Gamma) := d_0(\Gamma) \cup d_1(\Gamma)$ the set of vertices of $\Gamma$ and $E(\Gamma) := \Gamma \setminus V (\Gamma)$ the set of edges of $\Gamma$. For $e \in E(\Gamma)$ we call $d_0(e)$ and $d_1(e)$  the initial and
terminal vertices of  the edge $e$. 

A \emph{morphism} $\alpha:\Gamma\longrightarrow \Delta$ of profinite graphs is a continuous map with $\alpha d_i=d_i\alpha$ for $i=0,1$. If $\alpha$ is injective, the image  is called a subgraph of profinite graph $\Delta$, if $\alpha$ is surjective, then $\Delta$ is called a quotient graph of $\Gamma$.

 A profinite graph is called connected if its every finite quotient graph is connected. 

Let $\Gamma$ be  connected profinite graph. If $\Gamma=\varprojlim
	\Gamma_i$ is the inverse limit of the finite graphs $\Gamma_i$, then it induces the inverse
	system  $\{\pi_1(\Gamma_i)=\widehat \pi_1^{abs}(\Gamma_i)\}$ of the pro-$p$ completions of the abstract (usual) fundamental groups $\pi_1^{abs}(\Gamma_i)$. So 
	the pro-$p$ fundamental group $\pi_1(\Gamma)$ can be defined as $\pi_1(\Gamma)=\varprojlim_i
	\pi_1(\Gamma_i)$. 
	
\begin{definition}	If $\pi_1(\Gamma)=1$, then $\Gamma$ is called a {\em pro-$p$ tree}. 
\end{definition}

If $v,w\in V(\Gamma)$, the smallest pro-$p$ subtree of $\Gamma$ containing $\{v,w\}$ is called the $geodesic$ connecting $v$ and $w$, and is denoted $[v,w]$ (the definition is in pag. 83 of \cite{bib5}). 

\medskip
By definition a pro-$p$ group $G$ $acts$ on a profinite graph $\Gamma$ if we have a continuous action of $G$ on the profinite space $\Gamma$, such that $d_0$ and $d_1$ are $G$-maps. We shall denote by $G_m$ the stabilizer of $m\in\Gamma$ in $G$.  The action of a pro-$p$ group $G$ on a pro-$p$ tree $T$ is $irreducible$, if $T$ is the unique minimal $G$-invariant pro-$p$ subtree of $T$. The action is said to be $faithful$ if the Kernel $K$ of the action is trivial. If $G$ acts on $T$ irreducibly then the resulting action of $G/K$ on $T$ is faithful and irreducible.

\begin{definition} \rm
 Let $G$ be a pro-$p$ group, and $T$  a pro-$p$ tree on which $G$ acts continuously. For $g\in G$,
 \begin{itemize}
 \item $g$ is $elliptic$ if it fixes a vertex in $T$; if a subgroup of $G$ fixes a vertex we also call it elliptic.
 \item $g$ is $hyperbolic$ if it does not fixes a vertex in $T$. 
 \end{itemize}
 \end{definition}
 
  Observe that if $g\in G$ is hyperbolic, and $\langle g\rangle$ is the subgroup generated by $g$, then $\langle g\rangle$ acts freely on $T$. By \cite[Lemma 3.11]{bib5} there exists a unique nonempty minimal $\langle g\rangle$-invariant pro-$p$ subtree $D_g\subseteq T$. This holds also for any non-elliptic subgroup of $G$.\\
  
  Note that in the classical Bass-Serre theory a group $G$ acting on a tree without global fixed point splits as an amalgamated free product or an HNN-extension over the stabilizer of an edge. This is not true in the pro-p case in general but  is true for finitely generated case. As we shall often use this result we state it here.
  
  \begin{theorem} \label{splitting}\cite[Theorem 4.2]{bib7}  Let $G$ be a finitely generated
  pro-$p$ group acting  on a pro-$p$ tree $T$ without   global fixed points. Then
  $G$ splits non-trivially as a free amalgamated pro-$p$ product or
  pro-$p$ HNN-extension over some stabiliser of an edge of $T$.\end{theorem}

\subsection{Free pro-$p$ products with amalgamation}

	\begin{definition}[\cite{bib4}, Section 9.2]
		Let $G_1$ and $G_2$ be pro-$p$ groups and let $f_i: H \longrightarrow G_i$ $(i=1,2)$ be continuous monomorphisms of pro-$p$ groups. An amalgamated free pro-$p$ product of $G_1$ and $G_2$ with amalgamated subgroup $H$ is defined to be a pushout of $f_i$ $(i=1,2)$
$$		\xymatrix{	H \ar[r]^{f_1}\ar[d]^{f_2} & G_1 \ar[d]^{\varphi_1}\\
				G_2 \ar[r]^{\varphi_2} & G } $$

		in the category of pro-$p$ groups, i.e., a pro-$p$ group $G$ together with continuous homomorphisms $\varphi_i:G_i \longrightarrow G$ $(i=1,2)$ satisfying the following universal property: for any pair of continuous homomorphisms $\psi_i:G_i \longrightarrow K$ $(i=1,2)$ into a pro-$p$ group $K$ with $\psi_1f_1=\psi_2f_2$, there exists a unique continuous homomorphism $\psi: G \longrightarrow K$ such that the following diagram is commutative:
		
		$$\xymatrix{	H \ar[r]^{f_1} \ar[d]^{f_2} & G_1 \ar[d]^{\varphi_1} \ar[ddrr]^{\psi_1}\\
				G_2 \ar[r]^{\varphi_2} \ar[drrr]_{\psi_2} & G \ar[drr]^(.35){\psi} \\ &&&K  }$$
		
		An amalgamated free pro-$p$ product will be denoted by $G=G_1 \amalg_H G_2$.
	\end{definition}	
Note that an amalgamated free pro-$p$ product  can be also defined by presentation
$$G_1\amalg_H G_2=\langle G_1,G_2\mid rel(G_i), f_1(h)=f_2(h), h\in H, i=1,2\rangle.$$	
	
	Following the abstract notion, we can consider $H$ as a common subgroup of $G_1$ and $G_2$ and think of $f_1$ and $f_2$ as inclusions. However, unlike the abstract case where the canonical homomorphisms \[\varphi_i^{abs}:G_i \longrightarrow G_1 \star_{H} G_2 \] $(i=1,2)$ are always monomorphisms (cf. Theorem I.1 in \cite{bib2}), the corresponding maps in the category of pro-$p$ groups \[\varphi_i:G_i \longrightarrow G_1 \amalg_{H} G_2 \] $(i=1,2)$ are not always injective, i.e. it is not always  proper (in the terminology of  \cite{bib4}). However, we can make it proper by replacing $G_1, G_2, H$ with their images in $G$ (as explained in  \cite[Chapter 4]{bib5}).  If $G_1 \amalg_{H} G_2$ is proper we shall identify $G_1$, $G_2$ and $H$ with their images in $G$ and say that $G$ splits as the amalgamated free pro-$p$ product $G=G_1 \amalg_H G_2$.

	A free pro-$p$ product with cyclic amalgamation is always proper (see \cite{R}).
	Throughout the paper all free pro-$p$ products with amalgamation will be proper. 
	
	\medskip
	If $G=G_1 \amalg_H G_2$ and $H=G_1$  then $G=G_2$ and we call such splitting fictitious. All considered  splittings in the paper will be non-fictitious.

	\medskip
	We define the standard trees $S(G)$ on which $G=G_1\amalg_HG_2$ acts.
	
	\begin{itemize}
  \item Let $G=G_1\amalg_HG_2$. Then
 the vertex set is $V(S(G))= \displaystyle G/G_1\cup G/G_2$,
  the edge set is $E(S(G))= G/H$, and
  the initial and terminal vertices of an edge $gH$ are
  respectively  $gG_1$ and $gG_2$.
  
  By \cite[Theorem 4.1]{bib5} $S(G)$ is a pro-$p$ tree and the quotient graph $S(G)/G$ is an edge with two vertices.
	
	\subsection{Pro-$p$ HNN-extensions} \label{5}

	\begin{definition}[\cite{bib4}, Section 9.4]
		Let $H$ be a pro-$p$ group and let $f:A \longrightarrow B$ be a continuous isomorphism between closed subgroups $A$, $B$ of $H$. A pro-$p$ $HNN$-extension of $H$ with associated groups $A$, $B$ consists of a pro-$p$ group $G=HNN(H,A,t)$, an element $t \in G$ called the stable letter, and a continuous homomorphism $\varphi: H \longrightarrow G$ with $t(\varphi(a))t^{-1}=\varphi f(a)$ and satisfying the following universal property: for any pro-$p$ group $K$, any $k \in K$ and any continuous homomorphism $\psi: H \longrightarrow K$ satisfying $k(\psi(a))k^{-1}=\psi f(a)$ for all $a \in A$, there is a continuous homomorphism $\omega: G \longrightarrow K$ with $\omega(t)=k$ such that the diagram
		$$\xymatrix{G  \ar[dr]^{\omega} & \\
				H \ar[u]^{\varphi} \ar[r]^{\psi} & K }$$
		is commutative. 
	\end{definition}
	
A pro-$p$ HNN-extension $HNN((H,A,t)$ has the following presentation $$HNN(H,A,t)=\langle H, t\mid rel(H), a^t=f(a), a\in A\rangle.$$

	In contrast with the abstract situation, the canonical homomorphism $\varphi:H \longrightarrow G=HNN(H,A,t)$ is not always a monomorphism, i.e. it is not always  proper (in terminology of  \cite{bib4}). However, we can make it proper by replacing $H, A$ and $f(A)$ with their images in $G$ (as explained in  \cite[Chapter 4]{bib5}).   Throughout the paper all pro-$p$ $HNN$-extensions will be proper and in this case we shall idetify $H,A$ and $f(A)$ with their images in $G$ and say that $G$ splits as a pro-$p$ HNN-extension $G=HNN(H,A,t)$.

	\bigskip

We define a standard pro-$p$ tree on which HNN-extension acts.

    \item Let $G=HNN(G_1, H,t)$. Then
 the vertex set is $V(S(G))= \displaystyle G/G_1$,
  the edge set is $E(S(G))= G/H$, and
  the initial and terminal vertices of an edge $gH$ are
  respectively  $gG_1$ and $gtG_1$.
  \end{itemize}

By \cite[Theorem 4.1]{bib5} $S(G)$ is a pro-$p$ tree and the quotient graph $S/G$ is just a loop.

\subsection{ $\mathbb{Z}_{p}$-splittings}

  
 \bigskip
Let $G$ be a pro-$p$ group, $C_{1}$ and $C_{2}$ be subgroups of $G$ isomorphic to $\mathbb{Z}_{p}$,  the group of $p$-adic integers. A $\mathbb{Z}_{p}$-splitting of $G$ is a splitting as non-fictitious free pro-$p$ product with infinite cyclic amalgamation or as a proper pro-$p$ HNN-extension with infinite cyclic associated subgroup. Consider the following two $\mathbb{Z}_{p}$-splittings for $G$: 
\begin{itemize}
\item $G=A_{1}\amalg_{C_{1}} B_{1}$ or $G=HNN(A_{1},C_{1},t_{1}))$.
\item $G=A_{2}\amalg_{C_{2}} B_{2}$ or $G=HNN(A_{2},C_{2},t_{2}))$.
\end{itemize}
 
Let $T_1, T_2$ be the standard pro-$p$ trees  corresponding to the first $\mathbb{Z}_{p}$-splitting and the second $\mathbb{Z}_{p}$-splitting, respectively. Two given $\mathbb{Z}_{p}$-splittings are called:

\begin{itemize}

\item $Elliptic-Elliptic:$ if $c_{1}$ is elliptic in $T_{2}$ and $c_{2}$ is elliptic in $T_{1}$.

\item $Hyperbolic-hyperbolic:$ if $c_{1}$ is hyperbolic in $T_{2}$ and $c_{2}$ is hyperbolic in $T_{1}$.

\item $Hyperbolic-elliptic:$ if $c_{1}$ is hyperbolic in $T_{2}$ and $c_{2}$ is elliptic in $T_{1}$.
\end{itemize}

We shall see in Section \ref{sec4} that the last possibility in fact does not occur.

\subsection{Normalizer}

The following Proposition was proved in \cite[Proposition 8.1]{bib7}.
\begin{proposition} \label{teo6}
Let $p$ be a prime number and $C$ an infinite cyclic pro-$p$  group. Then:

\begin{itemize}

\item If the $\mathbb{Z}_{p}$-splitting is $G=A\amalg_{C} B$, then $N_G(C)=N_{A}(C)\amalg_{C} N_{B}(C)$.

\item If the $\mathbb{Z}_{p}$-splitting is $G=HNN(A,C,t)$, then:

\begin{itemize}

\item If $C$ and $C^{t}$ are  conjugate in $A$, then $N_{G}(C)=HNN(N_{A}(C),C,t')$ and $G=HNN(A,C,t')$.

\item If $C$ and $C^{t}$ are not conjugate in $A$ then $N_{G}(C)=N_{A^{t^{-1}}}(C)\amalg_{C} N_{A}(C)$.

\end{itemize}

\end{itemize}

\end{proposition}

\begin{proposition}\label{normalizer of cyclic}(\cite[Proposition 8.2]{bib7}) Let $G$ be a pro-$p$ group acting on a pro-$p$ tree $T$ and  $U$ be a cyclic
     subgroup of $G$ that does not stabilize any edge. Then one of the following happens:
     \begin{enumerate}
      
      \item[(1)] For some $g\in G$ and vertex $v$, $U\leq G_v$: then $N_G(U)=N_{G_v}(U)$.

\item[(2)] For all $g\in G$ and vertex $v$, $U\cap G_v=\{1\}$. Then $N_G(U)/K$ is
    either  isomorphic to $\Z_p$ or to a  dihedral pro-2 group $\Z/2\Z\amalg \Z/2\Z$, where $K$ is some normal subgroup of $N_G(U)$ contained in the stabilizer of an edge. 
        
       \end{enumerate}
     
     \end{proposition}

\section{Cyclic amalgamation}

In this section we shall consider  a $\Z_p$-splitting   $G=A\amalg_{C}B$ (resp. $G=HNN(A,C,t))$ as an amalgamated free pro-$p$ product or HNN-extension. Assuming that $N_G(H)$ is cyclic for every nontrivial $H\leq C$ we show 2-acylidricity of the action of $G$ on its standard pro-$p$ tree.

\begin{definition} Let $G$ be a pro-$p$ group and $H$  a subgroup of $G$. We say that $H$ is  \it{malnormal} in $G$ if  $H\cap H^{g}=1$, for any $g\in G-H$.
\end{definition}

\begin{proposition} \label{lemma help teo7}  Let $G$ be a pro-$p$ group and  $G=A\amalg_{C}B$ (resp. $G=HNN(A,C,t))$ be a $\mathbb{Z}_{p}$-splitting of $G$.
Suppose  $N_{G}(H)$ is  cyclic for every nontrivial open subgroup $H$ of $C$. Then $C$  is  malnormal  in $A$ or $B$ (resp. $C$ or $C^{t}$ is  malnormal  in $A$). 
\end{proposition}

\begin{proof} Case 1. $G=A\amalg_{C}B$. By cotradiction assume that $C$  is not malnormal in $A$ and $B$. Then there exist $a\in A-C$, $b\in B-C$ such that $C^{a}\cap C\neq 1\neq C\cap C^{b}$.  Consider $H=C\cap C^{a}\cap C^{b}$. Since  $[C:H]=[C:H^{a}]=[C: H^{b}]$ (as can be seen looking at finite quotients of $G$) one deduces that $a\in N_{A}(H)$, $b\in N_{B}(H)$.  By Proposition \ref{teo6} $N_{G}(H)=N_{A}(H)\amalg_{H} N_{B}(H)$.  Therefore $N_{G}(H)$ is not cyclic.\\
	
Case 2.	
Now suppose that $G=HNN(A,C,t)$. Assume on the contrary that  $C$ and $C^{t}$ are both not  malnormal  in $A$. Then there exist $a_1\in A-C$, $a_2\in A-C^{t}$  such that $C^{a_1}\cap C\neq 1\neq C^{ta_2}\cap C^{t} $. Put $H=C\cap C^{a_1}\cap C^{ta_2t^{-1}}\ $. 
  Since in every finite group conjugate subgroups have the same order, looking at finite quitients of $G$ one deduces that $[C:H]=[C:H^{a_1}]=[C^t:H^{ta_2}]$ and so $H=H^{a_1}=H^{ta_2t^{-1}}$, i.e. $a_1, a_2^{t^{-1}}\in N_{A}(H)$. By Proposition \ref{teo6} $N_{G}(H)$ can be a free amalgamated product $$N_{G}(H)=N_{A}(H)\amalg_{H}N_{A^{t^{-1}}}(H))$$ or a HNN-extesion $$N_{G}(H)=HNN(N_{A}(H),H,t).$$ In the first case observing that $a_2\in N_{A^{t^{-1}}}(H))$ we see that  $$N_{G}(H)=N_{A}(H)\amalg_{H}N_{A^{t^{-1}}}(H)$$  is not cyclic. In the second case $$N_{G}(H)=HNN(N_{A}(H),H,t)$$ is not cyclic because $N_{A}(H)$ is non-trivial. Thus we arrived at contradiction with the hypothesis and the proposition is proved.
\end{proof}

\begin{definition}\label{def: Profinite acylindrical}
The action of a pro-$p$ group $G$ on a pro-$p$ tree $T$ is said to be \emph{$k$-acylindrical}, for $k$ a constant, if the set of fixed points of $g$ has diameter at most $k$ whenever $g\neq 1$.
\end{definition}

\begin{corollary}\label{2-acy} Let $G$ be a pro-$p$ group such that $G=A\amalg_{C}B$ (resp. $G=HNN(A,C,t))$ be a $\mathbb{Z}_{p}$-splitting of $G$.
Suppose  $N_{G}(H)$ is cyclic for any nontrivial open subgroup $H$ of $C$. Then the action of $G$ on the standard pro-$p$ tree $S(G)$ of this splitting is 2-acylindrical.

\end{corollary}

\begin{proof} In case of an amalgamated free pro-$p$ product this is a direct consequence of Proposition \ref{lemma help teo7} and definitions of the standard pro-$p$ tree for a free amalgamated pro-$p$ product (cf. Sections 2.2). Indeed, if $g$ fixes three edges, it fixes three consecutive edges $e_1,e_2,e_3$ by \cite[Corollary 3.8]{bib5} and then conjugating $g$ if necessary we may assume that $e_2=1\cdot C$ (cf. the end of Sections 2.2) and so $e_1=aC$, $e_3=bC$ for some $a\in A,b\in B$. Then $g\in C\cap C^a\cap C^b=1$ by  Proposition \ref{lemma help teo7}.

Suppose  now $G$ is a pro-$p$ HNN-extensions $G=HNN(A,C,t)$. If $C$ and $C^t$ are conjugate in $A$ then by Proposition \ref{teo6} the normalizer  $N_G(C)\leq N_{G}(H)$ is not cyclic. Otherwise, $C\cap C^{at}=1$ for any $a\in A$ (since otherwise $t$ normalizes $C\cap C^{at}$ and so $N_G(C\cap C^{at} )$ is not cyclic). Again if $g\in G$ fixes three edges, it fixes three consecutive edges $e_1,e_2,e_3$ by \cite[Corollary 3.8]{bib5} and then conjugating $g$ if necessary we may assume that $e_2=1\cdot C$ (cf. the end of Sections 2.3) and so $e_1=a_1t^{-\epsilon}C$, $e_3=ta_2t^{-\epsilon}C$ for some $a_1,a_2\in A$, $\epsilon =0,1$ with $a_1t^{-\epsilon},ta_2t^{-\epsilon}\not\in C$.

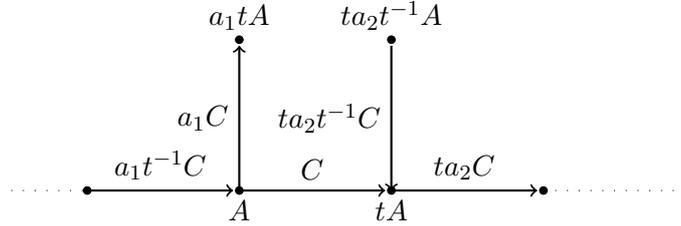
	$$\begin{tikzpicture}
	\node[circle,draw=black, fill=black, inner sep=0pt, minimum size=1.0mm] at (2,0) {};
	\node[circle,draw=black, fill=black, inner sep=0pt, minimum size=1.0mm] at (4,0) {};
	\node[circle,draw=black, fill=black, inner sep=0pt, minimum size=1.0mm] at (6,0) {};
	\node[circle,draw=black, fill=black, inner sep=0pt, minimum size=1.0mm] at (8,0) {};
	\node[circle,draw=black, fill=black, inner sep=0pt, minimum size=1.0mm] at (4,2) {};
	\node[circle,draw=black, fill=black, inner sep=0pt, minimum size=1.0mm] at (6,2) {};
	\draw (2,0)node[below] {};
	\draw (4,0)node[below] {$A$};
	\draw (6,0) node[below] {$tA$};
	\draw (4,2)node[above] {$a_1tA$};
	\draw (6,2)node[above] {$ta_2t^{-1}A$};
	\draw[thick][->] (2,0) -- node[above] {$a_1t^{-1}C$} (3.92,0);
	\draw[thick][->] (4,0) -- node[above] {$C$} (5.92,0);
	\draw[thick][->] (6,0) -- node[above] {$ta_2C$} (7.92,0);
	\draw[thick][->] (4,0) -- node[left] {$a_1C$} (4,1.92);
	\draw[thick][<-] (6,0) -- node[left] {$ta_2t^{-1}C$} (6,1.92);	
	\draw[loosely dotted] (1, 0) -- (2, 0);
	\draw[loosely dotted] (8, 0) -- (10, 0);
	\end{tikzpicture}$$

Since as was mentioned above $C\cap C^{at}=1$ for any $a\in A$ and  $ C\cap C^{a_1}\cap C^{ta_2t^{-1}}=1$ by  Proposition \ref{lemma help teo7},  $g$ has to be 1.
\end{proof}

\section{Excluding Hyperbolic-Elliptic $\mathbb{Z}_{p}$-Splitting }\label{sec4}


\bigskip
\begin{proposition}\label{trivial edge stabilizers} Let $G$ be a finitely generated pro-$p$ group admitting $\Z_{p}$-splitting $G=A\amalg_{C} B$ (resp. $G=HNN(A,C, t)$). If $A$ admit an action on an infinite pro-$p$ tree with trivial edge stabilizers where $C$ is elliptic, then so does $G$.

\end{proposition}

\begin{proof} By \cite[Theorem 9.6.1]{bib3} a finitely generated pro-$p$ group acts on an infinite pro-$p$ tree with trivial edge stabilizers if and only if it splits as a non-trivial free pro-$p$ product; we shall use it freely in this proof. 

\medskip
Suppose  $A$ acts on a pro-$p$ tree $T$ with trivial edge stabilizers. By \cite[Theorem 9.6.1]{bib3} $A$ is a non-trivial free pro-$p$ product $A=\coprod_{v\in V} A_v\amalg F$, where  $V$ is a transversal of $V(T)/G$ in $V(T)$ and $F$ a free pro-$p$ group acting freely on $T$. Since $C$ is elliptic in $T$, it is conjugate to $G_w$ for some $w\in V$, so w.l.o.g we may assume that $C\leq G_w$. 

\medskip
If $G=A\amalg_{C} B$, then we can rewrite $G$ as $(\coprod_{v\in V\setminus\{w\}} A_v\amalg F)\amalg (G_w\amalg_C B)$ and the proposition is proved in this case.

\bigskip
Suppose now $G=HNN(A,C, t)$. Then there exists $u\in V$ and $a \in A$ such that  $(C^t)^a$ is conjugate into $A_u$ and so replacing $t$ with $ta$ we may assume that $C^t\in A_u$. Then $G= (\coprod_{v\in V\setminus\{u,w\}} A_v\amalg F)\amalg HNN(A_w\amalg A_u, C,t)$. If $V\neq \{w,u\}$ then we are done, so we may assume that $G=HNN(A_w\amalg A_u, C,t)$ and so   $u\neq w$ since otherwise $\amalg_{v\in V(\Gamma)-\{u=w\}} \mathcal{A}(v)\neq 1$ because the splitting of  $A$ into the free product above is non-trivial. Hence $G=A_w\amalg_C A_u^{t^{-1}}\amalg \langle t\rangle $ (as follows from the presentation of HNN-extension) and so $G$  acts on a pro-$p$ tree with trivial edge stabilizers (cf. \cite[Section 4]{bib5}. The proposition is proved.

\end{proof}

\begin{theorem}\label{teo.gen.vpro-p}
   Let $G$ be a finitely generated  pro-$p$ group that does not split as a free pro-$p$ product. Then any two $\mathbb{Z}_{p}$-splittings of $G$ are either elliptic-elliptic or hyperbolic-hyperbolic.
\end{theorem}

\begin{proof}
    We are going to prove it by contradiction.  Let $G=A_i\amalg_{C_i} B_i$ or $G=HNN(A_i,C_i, t_i)$, $i=1, 2$ be two $\Z_p$-splittings and $T_1$, $T_2$ their standard pro-$p$ trees such that  w.l.o.g $c_1$ is hyperbolic in $T_2$ and $c_2$ is elliptic in $T_1$. 

 Since $c_2$ is elliptic in $T_1$, $C_2$ stabilizes a vertex of $T_1$. By Proposition \ref{trivial edge stabilizers} combined with Theorem \ref{splitting}  $A_2$ must fix a vertex $v$ in $T_1$.

 Now we will consider  2 possible cases.
    
\bigskip
    \textbf{Case 1:} (The second splittings is an amalgamation, say $G=A_2\amalg_{C_2}B_2$)\\
     By symmetry  $B_2 $ stabilize a vertex in $T_1$ as well. 
  Hence $A_2, B_2$ are contained in some conjugate of $A_1$ or $B_1$, say that $A_2\leq A_1$ (note that for HNN-extension case of the first splitting  there is no $B_1$). Then $B_2$ also has to be in $A_1$ (because otherwsie $A_2\cap B_2=1$ by \cite[Corollary 3.8]{bib5}). Hence $A_2$, $B_2\leq A_1$ and so $G=A_2\amalg_{C_2}B_2\leq A_1$, a contradiction.\\

\textbf{Case 2:} (The second splitting is a pro-$p$ HNN-extension)\\

Since  $A_2$ fixes a vertex in $T_1$,   w.l.o.g  we may assume that $A_2\leq A_1$, and we are going to prove that $t_2$ is in $A_1$. Since $c_{2}$ and $c^{t_2}_2$ are in $A_{2}<A_{1}$,   $c_2\in (A_1)^{t^{-1}_2}$. Hence $c_2\in A_1\cap A_1^{t^{-1}_2}$. If $t_2\not\in A_1$, by \cite[Theorem 4.3 (b),(c)]{bib5} $c_2\in A_1\cap A_1^{t^{-1}_2}<C^{a}_1$, for some $a\in A_1$. Then $c_2^{a^{-1}}\in C_1$, which is absurd, because $c_1$ is hyperbolic in $T_2$.

Thus $A_2$ and $t_2$ are in $A_1$, and $G=HNN(A_2,C_2,t_2)=\langle A_2,t_2\rangle \leq A_1$. So $G\leq A_{1}$, a final contradiction. \\
\end{proof}

\section{Hyperbolic-Hyperbolic $\mathbb{Z}_p$-splitting}\label{sec5}


Our objective  in this section is to prove the pro-p case of the Theorem 3.6 of \cite{bib1}.
Thus  in this section we fix two hyperbolic-hyperbolic $\mathbb{Z}_{p}$-splittings
 $G=A_{1}\amalg_{C_{1}}B_{1}$ ($G=HNN(A_{1},C_{1},t_{1})$) and $G=A_{2}\amalg_{C_{2}}B_{2}$ ($G=HNN(A_{2},C_{2},t_{2})$). 
As in the preceding section  $T_{1}, T_2$ stand for the standard pro-$p$ trees of the first and second  $\mathbb{Z}_{p}$-splittings, respectively.

\medskip
Note that if  $c_{1}$ is hyperbolic element of $T_{2}$  there exist a unique minimal $C_1$-invariant pro-$p$ $D_{1}$ of $T_{2}$, such that $C_{1}$ acts irreducibly on $D_{1}$ (see \cite[Lemma 3.11]{bib5}). Hence  $D_{1}$ is $N_{G}(C_{1})$-invariant. Since $C_{1}\triangleleft_{c}N_{G}(C_{1})$ and $C_{1}$ acts irreducibly on $D_{1}$, so does   $N_{G}(C_{1})$  (by Remark 4.2.1 (b) \cite{bib3}). Let $1\neq H_1\leq C_1$. Denoting by $K_{1}$ the kernel of the action of  $N_{G}(H_{1})$  on $D_{1}$ we deduce that  $N_{G}(H_{1})/K_{1}$ acts irreducibly and faithfully on $D_{1}$. Observe that as $D_{1}\subseteq T_{2}$, then $K_{1}\leq C^{g}_{2}$ for some $g\in G$, that is $K_{1}= 1$ or $K_{1}\cong \mathbb{Z}_{p}$. Moreover, replacing the second splitting by the $g$-conjugate we may assume for the rest of the section that $K_{1}\leq C_2$. Of course   similarly $C_2$ acts on a unique minimal $C_2$-invariant subtree $D_{2}$ of $T_1$ and all said above symmetrically holds here.

We split our consideration into two subsections $p>2$ and $p=2$.

\subsection{$p>2$}

In this subsection we consider $p>2$; in this case the idea of the proof is more explicit and proofs our more elegant.

The following Proposition is the pro-$p$ analog of the Proposition 3.3 of \cite{bib1}.

\begin{proposition} \label{prop1}  Let $H_i\neq 1$ be a subgroup of $C_i$. Then  $N_{G}(H_{i})$  has one of the following types:

\begin{enumerate}
\item[(i)]  cyclic group $\Z_p$; 
\item[(ii)]  
$\Z_p\times\Z_p$; in this case $C_1$ and $C_2$ commute.

\end{enumerate}
\end{proposition}

\begin{proof}
Put $N=N_{G}(H_{1})$ and let $K_1$ be the kernel of its action on $D_1$. By Proposition \ref{normalizer of cyclic}    $N/K_{1}\cong \mathbb{Z}_{p}$.
If  $K_{1}=1$ then we have item (i). 

\medskip Suppose
  $K_{1}\cong \mathbb{Z}_p$. Since $K_1$ is open normal in $C_2$ it acts feely and so irreducibly on $D_2$. Hence  $N_G(K_1)$ acts irreducibly on $D_2$. Note that $C_2,N_G(H_1)\leq N_G(K_{1})$ and $N_G(K_{1})$  can not act on $D_2$ faithfully by Proposition \ref{normalizer of cyclic} (2), so there exists a non-trivial kernel $K$ of this action and $N_G(K_{1})\cong K\rtimes \Z_p$. It follows that $N$ is open in $N_G(K_{1})$ and  $K\cap K_1=1$ since $K_1\leq C_2$ acts freely on $D_2$.
   Since $K_1$ is normal in $N$ it must centralize $K\cap N$.  For $p>2$ this means  $N_G(K_1) =K\times \Z_p$ (indeed,  a non-trivial action on $K\cong \Z_p$ is by multiplication by units and since $Aut(\Z_p)\cong \Z_p\times \Z/(p-1)\Z$ the action is faithful). Hence $N_G(K_1)$ is abelian; in particular $C_1$ and $C_2$ commute.  
\end{proof}

\begin{theorem}\label{teo7}
Let  $G$ be a finitely generated  pro-$p$ group that does not split as a free pro-$p$ product.  Let $G=A_{1}\amalg_{C_{1}} B_{1}$ (or $G=HNN(A_{1},C_{1},t_{1})$), and $G=A_{2}\amalg_{C_{2}} B_{2}$ (or $G=HNN(A_{2},C_{2},t_{2})$) be two hyperbolic-hyperbolic $\mathbb{Z}_{p}$-splittings of $G$. Suppose that $N_{G}(C_1)$ is not cyclic. Then $G\cong \Z_p\times \Z_p$.
\end{theorem}

\begin{proof}
We begin the proof with showing that the first $\Z_p$-splitting is an HNN-extension. Indeed, by Proposition \ref{prop1} $C_2\leq N_{G_1}(C_1)\cong \Z_p\times \Z_p$ and so 
we can write $N_{G}(C_i)$ as a pro-$p$ HNN-extension $HNN(\Z_p,\Z_p,t')\cong \Z_p\times \Z_p$ for $i=1, 2$, but  not as a non-fictiotious free amalgamated pro-$p$ product. Moreover, if the splitting of $N$ as a pro-$p$ amalgamated free product is fictiotious, i.e. either $N_{A_1}(C_1)=C_1$ or $N_{B_1}(C_1)=C_1$, say $N_{B_1}(C_1)=C_1$, then $N_G(C_1)=N_{A_1}(C_1)$ and so $C_2\leq A_1$ contradicting that $C_2$ is hyperbolic with rtespect to the first splitting. Hence by Proposition \ref{teo6} the first $\mathbb{Z}_p$-splittings of $G$ have to be an HNN-extension   and we can write  $G=HNN(A_1,C_1,t_1)$ such that $N_{G}(C_1)=HNN(N_{A_1}(C_1), C_1,t_1)=N_{A_1}(C_1)\amalg_{C_1}(C_1\times \langle t_1\rangle)$, i.e $C_1$ is normalized by $t_1$.  Then $$G=A_1\amalg_{C_1} (C_1\times \langle t_1\rangle)=A_1\amalg_{N_{A_1}(C_1)}N_{A_1}(C_1)\amalg_{C_1}(C_1\times \langle t_1\rangle)=A_1\amalg_{N_{A_1}(C_1)} N_{G}(C_{1}).$$ Thus $G$ can be rewritten  as  a free amalgamated pro-$p$ product 
\begin{equation}\label{eq1}
G=A_1\amalg_{N_{A_1}(C_1)} N_{G}(C_{1})
\end{equation} 

Moreover, since $C_2\leq N_G(C_1)$, $C_2$ is elliptic in $S_1$, where $S_1$ is the standard pro-$p$ tree of (\ref{eq1}). Since $C_1$ acts freely on $T_2$ and $A_1$ does not intersect any conjugate of $C_2$, by Proposition \ref{normalizer of cyclic} the normalizer $N_{A_1}(C_1)$ is  infinite cyclic. Thus if (\ref{eq1}) is non-fictitious, then (1) and the second splitting is a pair of  elliptic-hyperbolic $\Z_p$-splittings of $G$    contradicting  Theorem \ref{sela}.  Hence the splitting (\ref{eq1}) is fictitious, i.e. $A_1=N_{A_1}(C_1)$ and so  $G=N_G(C_1)$. Finally by Proposition \ref{prop1} (ii) $G\cong \Z_p\times \Z_p$.

\end{proof}

Combining this theorem with Proposition \ref{lemma help teo7} we deduce

 \begin{corollary}\label{malnormal p>2} With the hypothesis of Theortem \ref{teo7} suppose $G$ is non-abelian.
 Then $C_i$  is  malnormal  in either $A_i$ or $B_i$ (resp. $C_i$ or $C_i^{t_i}$ is  malnormal  in $A_i$), $i=1,2$. 
 
 \end{corollary}

\begin{theorem}\label{2-acylidrical}
Let   $G$ be a non-abelian finitely generated  pro-$p$ group that does not split as a free pro-$p$ product.  Let $G=A_{1}\amalg_{C_{1}} B_{1}$ (or $G=HNN(A_{1},C_{1},t_{1})$), and $G=A_{2}\amalg_{C_{2}} B_{2}$ (or $G=HNN(A_{2},C_{2},t_{2})$) be two hyperbolic-hyperbolic $\mathbb{Z}_{p}$-splittings of $G$.  Then the action of $G$ on the standard pro-$p$ trees of these splittings is 2-acylindrical.
\end{theorem}

\begin{proof} Suppose on the contrary the action of $G$ on say $T_1$ is not $2$-acylidrical.  Then by Proposition \ref{2-acy} $N_G(H)$ is not cyclic for some non-trivial subgroup $H$ of $C_1$. Hence by Proposition \ref{prop1} (ii) $N_G(C_1)$ is not cyclic. Therefore by  Theorem \ref{teo7}  $G$ is  abelian, contradicting the hypothesis. 

\end{proof} 

\subsection{$p=2$}

We continue with the settings of the beggining of this section.
The following Proposition is the pro-$2$ analog of  Propositions 3.1 and  3.3 of \cite{bib1}.

\begin{proposition} \label{prop1 p=2}   Let $H_i\neq 1$ be a subgroup of $C_i$. Then  $N_{G}(H_{i})$  has one of the following types:

\begin{enumerate}
\item[(i)]  cyclic group $\Z_2$; 
\item[(ii)] $\mathbb{Z}_{2}\rtimes \Z/2$ infinite dihedral pro-$2$ group;
\item[(iii)] $\Z_2\times \Z_2$;
\item[(iv)] 
$N_1\amalg_{K_1}N_2$, $[N_1:K_1]=2=[N_2:K_1]$. Moreover in this case either 
\begin{enumerate}
\item[(a)] $N_1\cong \Z_2\cong N_2$ and $N\cong\Z_2\rtimes\Z_2$ the pro-2 Klein bottle, or 
\item[(b)] $N_1\cong \Z/2\Z\amalg \Z/2\Z\cong N_2$, 
$N\cong (\Z_2\times \Z_2)\rtimes \Z/2\Z$, where $\Z/2\Z$ acts by inversion, or
\item[(c)]  $N_1\cong \Z/2\Z\amalg \Z/2\Z$, $N_2\cong\Z_2$ and  $N\cong (\Z_2\rtimes \Z_2)\rtimes \Z/2\Z$ is a semidirect product of the pro-$2$ Klein bottle and the group of order 2.
\end{enumerate}
\end{enumerate}
\end{proposition}

\begin{proof}
Consider $N=N_{G}(H_{1})$ and let $K_1$ be the kernel of its action on $D_1$. By Proposition \ref{normalizer of cyclic}   one of the following holds:
\begin{itemize}
\item $N/K_{1}\cong \mathbb{Z}_{2}$;
\item $N/K_{1}\cong \mathbb{Z}_{2}\rtimes \Z/2\Z$ infinite dihedral pro-$2$ group.
\end{itemize}
If  $K_{1}=1$ then we have item (i) or (ii). 

If $K_1\neq 1$ and 
  $N/K_{1}\cong \mathbb{Z}_2$ then we have  $N\cong K_1\rtimes \Z_2$. We show that if $N$ is non-abelian then  this is  a pro-2 Klein bottle. Indeed, if not then the action is faithful ($Aut(\Z_2)\cong \Z_2\times \Z/2\Z$) and so $N$ has a unique maximal cyclic normal   subgroup; since $K_1\leq C_2$ acts on $D_2$, so does $N$ but it can not act on $D_2$ faithfully by Proposition \ref{normalizer of cyclic} (2), so there exists non-trivial kernel $K$ of this action and $K\cap K_1=1$ since $K_1\leq C_2$ acts freely on $D_2$; this contradicts however the uniqueness of a maximal cyclic normal   subgroup of $N$. Thus either $N\cong\Z_2\times \Z_2$ or isomorphic to the pro-$2$ Klein bottle.  
 
 \medskip 
  Finally suppose that $K_1\neq 1$ and $N/K_{1}\cong \mathbb{Z}_{2}\rtimes \mathbb{Z}/2\mathbb{Z} $ is infinite dihedral. In this case $N=N_1\amalg_{K_1}N_2$ is a free amalgamated pro-2 product with $[N_1:K_1]=2=[N_2:K_1]$, i.e. (iv) holds. In this case if $N$ is torsion free, then $N_1\cong N_2\cong \Z_2$ and so it is the Klein bottle $\Z_2\rtimes \Z_2$, i.e.  we have  (iv)(a). 
  If both $N_1$ and $N_2$ are infinite dihedral pro-$2$ then (iv)(b) holds. If one of $N_i$, say $N_1$ is infinite dihedral pro-2 and $N_2$ is torsion free then the normal closure of $N_2$ in $N$ is the pro-$2$ Klein bottle, so (iv)(c) holds.
  
  \medskip 
  Note that $N_i$ for $i=1,2$ cannot be $\Z_2\times \Z/2\Z$. Indeed,  suppose without loss of generality that $N_1=K_1\times\Z/2\Z $. Since $K_1$ acts  freely on the pro-$p$ tree $D_2$,  $N_1$ acts  irreducibly on $D_2$. Hence by \cite[Theorem 3.12]{bib5} $\Z/2\Z$ is  the kernel of the action $K$ of $N_1$ on $D_2$. But  $K$ is contained in some conjugate of $C_1$ and so $K\cong \Z_2$, a contradiction. So $N_1\neq\Z_2\times \Z/2\Z$.

\end{proof}

\begin{remark}\label{virtually commute} The subcases (a),(b) and (c) of Proposition \ref{prop1 p=2} (iii) are  the pro-2 completion of the Klein bottle, a euclidean 4-branched sphere and an euclidean 2-branched (real) projective plane respectively. In particular in all this cases $N$ has a normal   subgroup of index 2 isomorphic to $\Z_2\times \Z_2$ and so $C_1$ and $C_2$ virtually commute.
\end{remark}

\begin{corollary} \label{metacyclic p=2} With the hypothesis  of Proposition \ref{prop1} assume $N_G(H_i)$ is torsion free. Then one of the following holds.
\begin{enumerate}
\item[(i)]  $N_{G}(H_{i})\cong \Z_2$;
\item[(ii)]  $N_{G}(H_{i})\cong \Z_2\times \Z_2$;
\item[(iii)] $N_{G}(H_{i})$ is isomorphic to the pro-2 Klein bottle $\Z_2\rtimes \Z_2$.
\end{enumerate}
\end{corollary}

\begin{lemma}\label{Lemma.pro-2} Let $G$ be a finitely generated pro-$2$ group and $G=A\amalg_{C} B$ ( resp. $G=HNN(A,C,t)$) is its $\Z_2$-splitting. Suppose $A$ admits an action on a pro-$2$ tree $T$ with edge stabilizers of order $\leq 2$ and without global fixed point, and $C$ is elliptic in $T$ (resp $C$ and $C^t$ are elliptic in $T_2$). Then $G$ splits over a  group of order $\leq 2$. (i.e, $G$ can be written as $G=G_1\amalg_{H} G_2$ or $HNN(G_1,H,t)$ with $|H|\leq 2$). 
	
\end{lemma}

\begin{proof}

	Since $G=A\amalg_{C}B$ (resp. $G=HNN(A,C,t)$) is finitely generated, the subgroup $A$ is finitely generated and $A$ acts on  a pro-$2$ tree $T$ without global fixed point having edge stabilizers of order at most 2 by hypothesis. Then by Theorem \ref{splitting} $A$ splits as an amalgamated free product $A=J_0\amalg_{H}J_1$ or as an HNN-extension $A=HNN(J_1,H,t)$ over a group $H$ of order at most 2. Moreover by \cite[Corollary 4.4]{bib7}  $C$ (resp. $C$ and $C^t$) is contained in $J_0$ or $J_1$ up to conjugation, say $J_1$, so w.l.o.g we may assume that $C\leq J_1$.  
	
\medskip	
	Now if $G=A\amalg_{C}B$, then  $G=J_0\amalg_{H}J_1\amalg_{C} B$ or $G=HNN(J_1,H,t_1)\amalg_C B=HNN(J_1\amalg_C B,H, t_1)$. 

\medskip
Suppose now	 $G=HNN(A,C,t)$. Since  $C^{ta}$ is contained in $J_0$ or $J_1$ for some $a\in A$, by replacing $t$ with $ta$ we may assume that $C^t\subseteq J_0\cup J_1$. Then $G$ has one of the following splittings:

1) If $A=J_0\amalg_{H}J_1$ and  $C, C^t\leq J_1$  then  $G=J_0\amalg_{H}HNN(J_1,C,t)$.

\smallskip
2) If $A=J_0\amalg_{H}J_1$ and  $ C^t\leq J_0$  then  $G=HNN(J_1\amalg_{C} J_0^{t^{-1}}), H, t^{-1})$.

\smallskip

3) If $A=HNN(J_1,H,t_1)$ then $G=HNN(HNN(J_1,C,t), H, t_1)$.
	 
	 The lemma is proved. 
	 
\end{proof}

\begin{lemma}\label{splitt} Let $G$ be a pro-$p$ and  $G=A_1\amalg_{C_1} B_1$ (or $G=HNN(A_1,C_1,t_1)$)  is a $\Z_p$- splitting of $G$. Then 

\begin{enumerate} 

\item[(i)] $G$ splits as an amalgamated free pro-$p$  product or as a HNN-extension over $N_{A_1}(C_1)$ in one of the following ways:
\begin{itemize}
	\item[a.] $G=A_1\amalg_{N_{A_1}(C_1)} (N_{G}(C_1)\amalg_{N_{B_1}(C_1)} B_1)$ if $G=A_1\amalg_{C_1} B_1$;
	\item[b.] $G=HNN(A_1\amalg_{N_{A_1}(C_1^{t_1})}N_{G}(C_1^{t_1}),N_{A_1}(C_1),t_1)$;
	\item[c.] $G=A_1\amalg_{N_{A_1}(C_1)}N_{G}(C_1)$, where $N_{G}(C_1)=HNN(N_{A_1}(C_1),C_1,t_1)$.
\end{itemize} 
\item[(ii)]  Let $G=A_{2}\amalg_{C_{2}} B_{2}$ (or $G=HNN(A_{2},C_{2},t_{2} )$)  be a  second splitting such that together with the first splitting it   is  a hyperbolic-hyperbolic pair.

\begin{itemize} 

\item[a.] Let $S_1$ be the standard pro-$p$ tree of any of the splittings in $(i)$ over $N_{A_1}(C_1)$. If $A_2$ and $B_2$ (resp. $A_2$) stabilize a vertex in $S_1$, then so is $G$.

\item[b.] If   $N_{A_1}(C_1)\neq A_1$, then $A_2$ or $B_2$ (resp. $A_2$) is not elliptic in the splittigs (i) over $N_{A_1}(C_1)$. Moreover, if in addition $N_{G}(C_1)$ is neither cyclic nor dihedral and $N_{A_1}(C_1)\neq A_1$ then $G$ splits over a group of order $\leq 2$.
\end{itemize}
\end{enumerate}

\end{lemma}

\begin{proof} $(i)$ Suppose the $\Z_p$-splitting is  an amalgamated free pro-$p$ product $$G=A_1\amalg_{C_1}B_1.$$ Then by Proposition \ref{teo6} $N_G(C_1)=N_{A_1}(C_1 )\amalg_{C_1} N_{B_1}(C_1)$ and thus $G$ admits a decomposition as follows:

\begin{equation}\label{eq2k}
	\begin{array}{cl} 
	G       &=A_1\amalg_{C_1}B_1\\
	&=\boxed{A_1\amalg_{N_{A_1}(C_1)} (N_{G}(C_1)\amalg_{N_{B_1}(C_1)} B_1)}.\\
	\end{array}
	\end{equation}
Thus $G$ splitts as a  free pro-$p$ product with $N_{A_1}(C_1)$ amalgamated  in this case.\\
	
On the other hand if $G=HNN(A_1,C_1,t_1)$ is an HNN-extension, by Proposition \ref{teo6} the normalizer $N_G(C_1)$ is an amalgamated free pro-$p$ product $N_G(C_1)=N_{A_1}(C_1)\amalg_{C_1} N_{A_1^{t_1^{-1}}}(C_1)$ or  an HNN-extension $N_G(C_1)=HNN(N_{A_1}(C_1),C_1,t_1)$. In the first case putting $C'=C_1^{t_1}$ we have that $N_{G}(C')=N_{A^{t_1}_1}(C')\amalg_{C'}N_{A_1}(C')$ and so

\begin{equation}\label{eq3k}
\boxed{G=HNN(A_1\amalg_{N_{A_1}(C')}N_{G}(C'),N_{A_1}(C_1),t_1).}.
\end{equation}

Indeed, since $N_{G}(C')$ is a subgroup of $G$ and $N_{A^{t_1}_1}(C')=(N_{A_1}(C_1))^{t_1}$ if and only if $C'=(C_1)^{t_1}$, we can eliminate $N_{G}(C_1)$ from the presentation \eqref{eq3k} to arrive at $\langle A_1, t_1: Rel(A_1),\, C'=C^{t_1}_1\rangle=HNN(A_1,C_1,t_1)$. This finishes the prove in this case.

In the second case 
\begin{equation}\label{eq4k}
	\begin{array}{cl} 
	G       &=HNN(A_1,C_1,t_1)\\
	&=\overline{\langle A_1,t_1: C_1^{t_1}=C_1\rangle}\\
	&=A_1\amalg_{C_1} C_1\times \overline{\langle t_1\rangle}\\
	&=A_1\amalg_{N_{A_1}(C_1)}N_{A_1}(C_1)\amalg_{C_1}C_1\times \overline{\langle t_1\rangle}\\
	&= A_1\amalg_{N_{A_1}(C_1)} HNN(N_{A_1}(C_1),C_1,t_1)\\
	&=\boxed{A_1\amalg_{N_{A_1}(C_1)}N_{G}(C_1)}
	\end{array}
	\end{equation}
This finished the proof of (i).\\

$(ii)$

(a.) Assume that the second $\Z_2$-decomposition is an amalgamated free product pro-$p$, $G=A_2\amalg_{C_2}B_2$. We are going to prove that if $A_2$ and $B_2$ stabilize a vertex in $S_1$, then $A_2$ and $B_2$ stabilize the same vertex. Indeed, if $A_2$ stabilizes $v\in V(S_1)$, $B_2$ must be in the stabilizer of $v$, since otherwise  $C_2=A_2\cap B_2$ is contained in a conjugate of $N_{A_1}( C_1)\leq A_1$ by \cite[Corollary 4.1.6]{bib3} contradicting the hypothesis.

\noindent On the other hand if the second $\Z_2$-decomposition is pro-2 HNN-extension $$G=HNN(A_2,C_2,t_2),$$ then $ A_2$ and $(A_2)^{t^{-1}_2}$ act on $S_1$. If $A_2$ stabilizes the vertex $v\in V(S_1)$, then $(A_2)^{t^{-1}_2}$, stabilizes a vertex in $S_1$. Then the element $t_2$ must  stabilize  $v$, since otherwise  by \cite[Corollary 4.1.6]{bib3} $C_2=A_2\cap (A_2)^{(t_2)^{-1}}$ is in a conjugate of $N_{A_1}(C_1)\leq A_1$ which is absurd.\\

\medskip

$b.$ As $A_1\neq N_{A_1}(C_1)$, then by (ii)a. the case (i)b. can not happen. The second factor of the decompositions (\ref{eq2k}) or (\ref{eq4k}) is a proper subgroup of $G$.\\

Assume that the second $\Z_2$-decomposition is an amalgamated free pro-$p$ product  $G=A_2\amalg_{C_2}B_2$. Note that $A_2$, and $B_2$ are not contained in a conjugate of $A_1$ because if it were, $C_2$ would be contained in a conjugate of $A_1$ contradicitng the hypothesis. On the other hand, if $A_2$ is contained in a conjugate of the second factor, then $B_2$ must be in the same conjugate that $A_2$ (by $(ii)a.$). Then $G=A_2\amalg_{C_2}B_2$ would be contained in some conjugate of the second factor  which is absurd because the second factor is a proper subgroup in $G$. Thus $A_2$ is not elliptic in $S_1$.\\

On the other hand suppose the second $\Z_2$-decomposition is a pro-2 HNN-extension, $G=HNN(A_2,C_2,t_2)$. We know that $A_2$ is not contained in a conjugate of $A_1$ because if it were, then $C_2$ would be contained in a conjugate of $A_1$ contradicting the hypothesis. If $A_2$ is contained in a conjugate of the second factor, then  $t_2$ must be in the same conjugate as $A_2$ (by $(ii)a.$). Hence $G=HNN(A_2,C_2,t_2)$ is contained in some conjugate of the second factor  which is absurd because the second factor is proper subgroup of $G$. Thus $A_2$ is not elliptic in $S_1$.	
	
	\medskip
	To prove the second statement of b. assume that  $N_G(C_1)$ is neither cylic nor dihedral. Then the kernel $K_1$ of its action on $D_1$ is non-trivial. The group $K_1$ is elliptic in $S_1$ and therefore so is $C_2$. Note that by Proposition \ref{normalizer of cyclic} (2) $N_{A_1}(C_1)$ is either cyclic or infinite dihedral. If the normalizer $N_{A_1}(C_1)$ is infinite cyclic, then by Theorem \ref{teo.gen.vpro-p} $G$ admits a decomposition as a free pro-$p$ product. If the normalizer $N_{A_1}(C_1)\cong \Z_2\rtimes \Z/2\Z$, then $N_{A_1}(C_1)$ contains a cyclic normal subgroup $C$ of index at most $2$ generated by a hyperbolic element in $T_2$. Then  $[N_{A_1}(C_1): C]\leq 2$  and so   considering the action of $A_2$ on $S_1$ we see that   $|A_2\cap N_{A_1}(C_1)|\leq 2$ since $A_2\cap C=1$. It follows that   $A_2$-stabilizers of edges of $S_1$ are either trivial or groups of  order 2. Since $A_2$ is not elliptic,     the group $G$ splitts over  a group of order $\leq 2$ by Lemma \ref {Lemma.pro-2}, as required.

\end{proof}

\begin{theorem}\label{teo8}
	Let be $G$ a finitely generated pro-$2$ group  which does not split over a cyclic group of order $\leq 2$.  Let $G=A_{1}\amalg_{C_{1}} B_{1}$ (or $G=HNN(A_{1},C_{1},t_{1})$), and $G=A_{2}\amalg_{C_{2}} B_{2}$ (or $G=HNN(A_{2},C_{2},t_{2})$) be two hyperbolic-hyperbolic $\mathbb{Z}_{2}$-splittings of $G$. Suppose that $N_{G}(C_1)$ is neither cyclic nor dihedral. Then  $G=N_G(C_1)=N_G(C_2)$ is virtually abelian isomorphic to one of  pro-2 groups listed in (iii) or (iv) of Proposition \ref{prop1 p=2}. 
	
\end{theorem}

\begin{proof}  By Lemma \ref{splitt} $(i)$ $G$ admits the decompositions \eqref{eq2k}-\eqref{eq4k} and since   by  Lemma  \ref{splitt} (ii)b.  $A_1=N_{A_1}(C_1)$ the decompositions reduce to the following form: 
		\begin{equation}\label{eqj2}
		G= (N_{G}(C_{1})\amalg_{N_{B_{1}}(C_{1})} B_{1})
		\end{equation}

 \begin{equation}\label{eqj3}
G=A_1\rtimes \langle t_1\rangle\amalg_{N_{A_1}(C_1^{t_1})}N_{G}(C_1^{t_1})\end{equation}

\begin{equation}\label{eqj4}
G=N_G(C_1)
\end{equation}

		 Suppose now the first $\mathbb{Z}_2$-splitting of $G$ is an  amalgamated free pro-$2$ product $G=A_1\amalg_{C_1} B_1$.  Then  swapping $A_1$ and $B_1$ in the first splitting of $G$  we also have  $N_{B_{1}}(C_{1})=B_1$. Thus $G=N_{G}(C_1)=N_{A_1}(C_1)\amalg_{C_1}N_{B_1}(C_1)$ is isomorphic  to one of the pro-$2$ groups  of Proposition \ref{prop1 p=2} $(iv)$. 
		
		\bigskip

		 If the first $\mathbb{Z}_2$-splitting of $G$ is a pro-$2$ HNN-extension, then  in the case of $N_G(C_1)=HNN(N_{A_1}(C_1), C_1, t)$  we have (\ref{eqj4}) by Lemma \ref{splitt}(i)c.
		 
\medskip		 
		 Otherwise we have  (\ref{eqj3}). Conjugating this equation by $t_1^{-1}$ we get 
		 
\begin{equation}\label{eqj5}
G=A_1\rtimes \langle t_1\rangle\amalg_{N_{A_1}(C_1)}N_{G}(C_1)
\end{equation}		 
		 
		  and  since $A_1=N_{A_1}(C_1)$ and $N_G(C_1)=N_{A_1}(C_1)\amalg_{C_1} N_{A_1^{t^{-1}}}(C_1)$ (cf. Lemma \ref{splitting}(i)), it rewrites as  
\begin{equation}\label{eqj6}
G=A_1\rtimes \langle t_1\rangle\amalg_{C_1}N_{A_1^{t^{-1}}}(C_1)
\end{equation}			  

This splitting forms hyperbolic-hyperbolic splitting with the second splitting as well, so applying to it the second paragraph of the proof we deduce that $G=N_G(C_1)$.	
	
	\medskip
	Thus, $G=N_{G}(C_1)$ and  we have either Case (iii) or  Case (iv)  of Proposition \ref{prop1 p=2}.

\end{proof}

\end{document}